\documentclass[11pt]{article}
\usepackage{graphicx}
\usepackage{amsmath}
\usepackage{amssymb}
\usepackage{epstopdf}

\def\c{\centerline}

\def\L{{\bf L}}

\def\D{{\mathcal D}}

\def\S{{\mathcal S}}

\def\Hat{\widehat}

\def\C{{\mathcal C}}
\def\U{{\mathcal U}}

\def\ds{\displaystyle}
\def\sqr#1#2{\vbox{\hrule height .#2pt
\hbox{\vrule width .#2pt height #1pt \kern #1pt
\vrule width .#2pt}\hrule height .#2pt }}
\def\square{\sqr74}
\def\endproof{\hphantom{MM}\hfill\llap{$\square$}\goodbreak}

\def\bega{\begin{array}}
\def\enda{\end{array}}
\def\begi{\begin{itemize}}
\def\endi{\end{itemize}}
\def\T{{\mathbb T}}
\def\O{{\cal O}}

\def\R{{\mathbb R}}
\def\Z{{\mathbb Z}}
\def\ov{\overline}

\def\v{\vskip 1em}
\def\vs{\vskip 2em}

\def\be{\begin{equation}}
\def\beq{\begin{equation}}
\def\bel{\begin{equation}\label}
\def\eeq{\end{equation}}

\begin{document}
\title{\bf Structurally Stable Singularities 
for a Nonlinear   Wave Equation}
\vs
\author{Alberto Bressan, Tao Huang, and Fang Yu\\ \\ 
Department of Mathematics, Penn State University\\  University Park,
Pa. 16802, U.S.A.\\ \\ e-mails: 
bressan@math.psu.edu$\,$,~txh35@psu.edu$\,$,~fuy3@psu.edu}
\maketitle

\c{\it Dedicated to Tai Ping Liu in the occasion of his 70-th birthday}
\v

\begin{abstract} 
For the nonlinear wave equation $u_{tt} - c(u)\big(c(u) u_x\big)_x~=~0$,
it is well known that solutions can develop singularities in finite time. For an open dense set of initial data, 
the present paper provides a detailed asymptotic
description of the solution in a neighborhood of each singular point,
where $|u_x|\to\infty$. The different structure of
conservative and dissipative solutions is analyzed.
\end{abstract}

\section{Introduction} 
\setcounter{equation}{0}
The nonlinear wave equation
\bel{1.1}
u_{tt} - c(u)\big(c(u) u_x\big)_x~=~0\,,\eeq
provides a mathematical model for the behavior of
nematic liquid crystals.
Solutions have been studied
by several authors \cite{BC1, BCZ, BH, BZ,  GHZ, HR,
ZZ1, ZZ2}.  We recall that, even for smooth initial data
\bel{1.2}
u(x,0)~=~u_0(x)\,,\qquad\qquad
u_t(x,0)~=~u_1(x)\,,
\eeq
regularity can be lost in finite time.
More precisely, the $H^1$ norm of the solution $u(\cdot,t)$
remains bounded, hence $u$ is always H\"older continuous, 
but the norm of the gradient $\|u_x(\cdot,t)\|_{\L^\infty}$ 
can blow up in finite time.

The paper \cite{BZ} introduced a nonlinear 
transformation
of variables that reduces (\ref{1.1}) to a semilinear system.
In essence, it was shown that
the quantities
$$w~\doteq~2\arctan\bigl(u_t +c(u) u_x\bigr),
\qquad\qquad z ~\doteq~2\arctan\bigl(u_t -c(u) u_x\bigr),$$
satisfy a first order semilinear system of equations, w.r.t.~new 
independent variables $X$, $Y$ constant along characteristics. 
Going back to the original  variables $x,t,u$, one obtains
a global solution of the 
wave equation (\ref{1.1}). 

Based on this representation and using ideas from 
\cite{D,DD,GG,G}, in \cite{BC1} it was recently 
proved that,
for generic initial data, 
the conservative solution is smooth outside a finite 
number of points and curves in the $t$-$x$ plane.   
Moreover, conditions were identified which guarantee the 
{\it structural stability} of the set of singularities.  
Namely, when these generic 
conditions hold, the topological structure of the singular set
is not affected by a small $\C^3$ perturbation of the initial data.

Aim of the present paper is to derive a detailed asymptotic
description of these structurally stable solutions, 
in a neighborhood of each singular point.
This is achieved both for conservative and for dissipative solutions
of (\ref{1.1}).
We recall that conservative solutions satisfy an additional conservation 
law for the energy, so that the total energy
$${\mathcal E}(t)~=~{1\over 2} \int [u_t^2 +  c^2(u) u_x^2]\, dx$$
coincides with a constant for a.e.~time $t$.
On the other hand, for dissipative solutions the total energy 
is a monotone decreasing function of time. 
A representation of dissipative solutions in terms of a suitable semilinear system in characteristic coordinates can be found in 
\cite{BH}.

The remainder of this paper is organized as follows.  In Section~2
we review the variable transformations introduced in \cite{BZ}
and the conditions for structural stability derived in \cite{BC1}.
Section~3 is concerned with conservative solutions.
In this case, for smooth initial data the map
\bel{map}
(X,Y) ~\mapsto ~(x,t,u,w,z)(X,Y)\eeq 
remains globally smooth, on the entire
$X$-$Y$ plane.  
To recover the singularities
of the solution $u(x,t)$ of (\ref{1.1}), it suffices to 
study the Taylor approximation of (\ref{map}) at points where
$w=\pi$ or $z=\pi$.  In Section~4 we perform a similar analysis
in the case of dissipative solutions.   This case
is technically more difficult,
because the corresponding semilinear system
has discontinuous source terms.  

We remark that, for conservative solutions, 
a general uniqueness theorem
has been recently established in \cite{BCZ}. 
On the other hand, for dissipative solutions no general 
result on uniqueness or continuous dependence is yet known.
Whether structurally stable dissipative solutions are {\it generic},
arising from an open dense set of $\C^3$ initial data, is 
also an open problem.

\section{Review of the equations} 
\setcounter{equation}{0}
Throughout the following, on the wave speed $c$ we assume
\begi
\item[{\bf (A)}] The map
 $c:\R\mapsto \R_+$ is
 smooth and  uniformly positive.   The quotient   $c'(u)/c(u)$ is uniformly 
 bounded. Moreover,
 the following generic condition is satisfied:
 \bel{morse}  c'(u)~=~0\qquad\implies\qquad c''(u)~\not=~ 0.\eeq 
 \endi
 Because of (\ref{morse}), the derivative $c'(u)$ 
can vanish only at isolated points.

In (\ref{1.2})   we consider initial data $(u_0, u_1)$ in the product space $H^1(\R)\times \L^2(\R)$.
It is convenient to introduce the variables
\beq
\left\{
\begin{array}{rcl}
R & \doteq  &u_t+c(u)u_x\,, \\
S & \doteq  &u_t-c(u)u_x\,,
\end{array} \right.\label{2.1}
\eeq
In a smooth solution, $R^2$ and $S^2$ satisfy the  balance laws
\bel{2.4}\left\{
\begin{array}{rcl}
(R^2)_t - (cR^2)_x & = & {c'\over 2c}(R^2S - RS^2)\, , \\ [3mm]
(S^2)_t + (cS^2)_x & = &  {c'\over 2c}(S^2R-SR^2 )\,.
\end{array}
\right.
\eeq
As a consequence,  the energy  is conserved:
\beq
E~\doteq ~{1\over 2}\big(u_t^2+c^2u_x^2\big)~=~{R^2+S^2\over 4}\,,
\eeq
One can think of $R^2$ and $S^2$ as the energy of backward and 
forward moving waves, respectively.  Notice that these are not separately conserved. 
Indeed, by (\ref{2.4}) energy can be exchanged between forward and backward waves.

A major difficulty in the analysis of (\ref{1.1}) is the possible
breakdown of regularity of solutions.   Indeed, even for smooth
initial data, the quantities $u_x, u_t$ can blow up in finite time.
To deal with possibly unbounded values of $R,S$, following \cite{BZ}
we introduce a new set of dependent variables:
\bel{wzdef}
w~\doteq~ 2\arctan R\,,\qquad\qquad z~\doteq~ 2\arctan S\,.
\eeq

\begin{figure}[htb]
\centering
\includegraphics[width=0.55\textwidth]{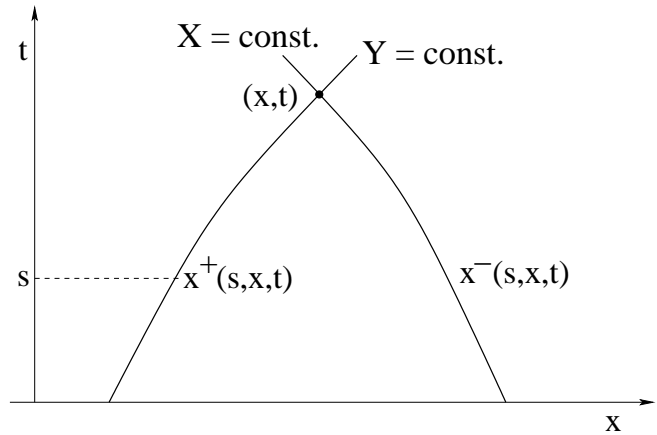} 
\caption{\small The backward and forward characteristic 
through the point $(x,t)$.}
\label{f:wa71}
\end{figure}

To reduce the equation (\ref{1.1}) to a semilinear one, 
it is convenient to
perform a further change of independent variables (Fig.~\ref{f:wa71}).
Consider
the equations for the forward and backward characteristics:
\beq
\dot x^+=c(u)\,,\qquad\qquad \dot x^-=-c(u)\,.\label{2.12} 
\eeq
The characteristics passing through the point $(x,t)$
will be denoted by
$$
s\mapsto x^+(s,x,t)\,,\qquad\qquad s\mapsto x^-(s,x,t)\,,
$$
respectively.
As coordinates $(X,Y)$ of a point $(x,t)$ we shall use the quantities
\beq
X~\doteq ~x^-(0,x,t)\,,\qquad\qquad Y~\doteq ~-x^+(0,x,t)
\,.\label{XY}
\eeq
For future use, we now introduce the further variables
\bel{pqdef}
p\doteq {1+R^2\over X_x}\,,\qquad\qquad q\doteq {1+S^2\over -Y_x}\,.
\eeq
Starting with the nonlinear equation (\ref{1.1}),
using $X,Y$ as independent variables one obtains a semilinear
hyperbolic system with smooth coefficients for the variables 
$u,w,z,p,q,x,t$, namely
\bel{u}
\left\{
\begin{array}{ccc}
u_X &= & {\sin w\over 4c} \, p\,,\\[4mm]
 u_Y &= & {\sin z\over 4c} \, q\,,
\end{array}\right. 
\eeq

\bel{wz}
\left\{
\begin{array}{ccc}
w_Y&=&{c'\over 8c^2}\,( \cos z - \cos w)\,q\,,\\[4mm]
z_X&=&{c'\over 8c^2}\,( \cos w - \cos z)\,p\,,
\end{array}\right. 
\eeq

\bel{pq}
\left\{
\begin{array}{ccc}
p_Y &= & {c'\over 8c^2}\, (\sin z-\sin w)\,pq\,,\\[4mm]
q_X &= & {c'\over 8c^2}\, (\sin w-\sin z)\,pq\,,
\end{array}\right.
\eeq

\bel{x}
\left\{
\begin{array}{ccr}
x_X&=&{(1+\cos w)\,p\over 4}\,,\\[4mm]
x_Y&=&-{(1+\cos z)\,q\over 4}\,,
\end{array}\right. 
\eeq

\bel{t}
\left\{
\begin{array}{ccr}
 t_X&=&{(1+\cos w)\,p\over 4c}\,,\\[4mm]
t_Y&=&{(1+\cos z)\,q\over 4c}\,.
\end{array}\right. 
\eeq
See \cite{BZ} for detailed computations.
Boundary data can be assigned on the  line 
$\gamma_0=\{(X,Y)\,;~ X+Y=0\}$,
by setting
 \bel{bdata}  \left\{
\begin{array}{rl}
u(s,\,-s)&=~\ov u(s)\,,\\
x(s,\,-s)&=~\ov x(s)\,,\\
t(s,\,-s)&=~\ov t(s)\,,\enda
\right.
\qquad\qquad
\left\{
\begin{array}{rl}
w(s,\,-s)&= ~\ov w(s)\,,\\ z(s,\,-s)&= ~ \ov z(s)\,,
\end{array}\right.\qquad \qquad
\left\{
\begin{array}{rcl}
p(s,\,-s)&=~ \ov p(s)\,,\\
q(s,\,-s)&=~\ov q(s)\,,
\end{array}\right.  
\eeq
for suitable smooth functions $\ov u,,
\ov x,\ov t, \ov w,\ov z,\ov p,\ov q$. 

\v
{\bf Remark 1.}  The above system is clearly 
invariant w.r.t.~the addition
of an integer multiple of $2\pi$ to the variables $w,z$.
Taking advantage of this property, in the following we 
shall regard $w,z$ as points in the quotient manifold
$\T\doteq \R/2\pi\Z$.  As a consequence, we have the implications
\bel{impl}\bega{rl}
w&\not=~\pi\qquad\implies\qquad \cos w~>~-1\,,\\[3mm]
z&\not=~\pi\qquad\implies\qquad \cos z~>~-1\,.\enda\eeq
\v
{\bf Remark 2.}  The system (\ref{u})--(\ref{t}) is overdetermined.
Indeed, the functions $u,x,t$ can be computed by
using either one of the equations in (\ref{u}), (\ref{t}), 
(\ref{x}), respectively. As shown in 
\cite{BC1}, in order that all the above equations be 
simultaneously satisfied along the line
$\gamma_0$ one needs the additional compatibility conditions
\bel{cc1} {d\over ds}  \ov u(s)~=~{\sin \ov w(s)\over 4c(\ov u(s))} \,
\ov p(s)- {\sin \ov z(s)\over 4c(\ov u(s))}\,\ov q(s)\,,\eeq
\bel{cc2}{d\over ds} \ov x(s) ~=~{(1+\cos \ov q(s)) \ov p(s) 
+ (1+\cos \ov z(s) )\ov q(s)\over 4}\,,\eeq
\bel{cc3}
{d\over ds} \ov t(s) ~=~{(1+\cos \ov  w(s)) \ov p(s) 
- (1+\cos \ov z(s) )\ov q(s)
\over 4 c(\ov u(s))}\,.\eeq
In turn, if (\ref{cc1})--(\ref{cc3}) hold along $\gamma_0$, then 
a unique solution to the  system (\ref{u})--(\ref{t}) 
can be constructed, on the entire $X$-$Y$ plane.
\v
Given initial data  $(u_0,u_1)$ in (\ref{1.2}), we assign
 boundary data  (\ref{bdata}) on the line $\gamma_0$, 
by setting
\bel{data}  \left\{
\begin{array}{rl}
\ov u(x)&=~u_0(x)\,,\\
\ov t(x)&=~0\,,\\
\ov x(x)&=~x\,,\enda
\right.
\qquad\quad
\left\{
\begin{array}{rl}
\ov w(x)&= ~2\arctan R(x,0)\,,\\ \ov z(x)&= ~ 2\arctan S(x,0)\,,
\end{array}\right.\quad \qquad
\left\{
\begin{array}{rcl}
\ov p(x) &\equiv & 1+ R^2(x,0)\,,\\
\ov q(x) &\equiv & 1+S^2(x,0)\,.
\end{array}\right.
\eeq
We recall that, at time $t=0$, by (\ref{1.2}) one has
$$\bega{lr}R(x,0) &=~(u_t + c(u) u_x)(x,0) ~=~
u_1(x) + c(u_0(x)) u_{0,x}(x),\\[4mm]
S(x,0) &=~(u_t - c(u) u_x)(x,0) ~=~u_1(x) - 
c(u_0(x)) u_{0,x}(x).\enda$$
As proved in \cite{BC1}, 
for any choice of $u_0,u_1$ in (\ref{data}) the compatibility
conditions (\ref{cc1})--(\ref{cc3}) are automatically satisfied.

The following theorems summarize the main results on 
conservative solutions, proved in \cite{BZ, BC1, BCZ}. 
As before, $\U$ denotes the product space in (\ref{U}).
\v
{\bf Theorem 1.} {\it Let the wave speed $c(\cdot)$ 
satisfy the assumptions
{\bf (A)}.  

 Given initial data $(u_0,u_1)\in H^1(\R)\times\L^2(\R)$, there exists
a unique solution\\ $(X,Y)\mapsto 
(u,w,z,p,q,x,t)(X,Y)$ to the system (\ref{u})--(\ref{t}) with 
boundary data (\ref{bdata}), (\ref{data}) 
assigned along the line $\gamma_0$.   Moreover, the set 
\bel{graph}
\hbox{\rm Graph}(u)~\doteq~\Big\{\bigl(x(X,Y), ~t(X,Y),~u(X,Y)\bigr)\,;
~~(X,Y)\in\R^2\Big\}\eeq
is  the graph of the unique conservative solution $u=u(x,t)$
of the Cauchy problem (\ref{1.1})-(\ref{1.2}).
}

\v
{\bf Theorem 2.} {\it  Let the assumptions
{\bf (A)} be satisfied and let $T>0$ be given.  
Then there exists
an open dense set 
\bel{U} \D~\subset~
\U~\doteq~ \Big(\C^3(\R)\cap H^1(\R)\Big) 
\times\Big(\C^2(\R)\cap\L^2(\R)\Big)\eeq
 such that the following holds.

For every initial data $(u_0,u_1)\in \D$, 
the corresponding solution $(u,w,z,p,q,z,t)$ 
of  (\ref{u})--(\ref{t}) with 
boundary data (\ref{bdata}), (\ref{data}) has level sets
$\{w=\pi\}$ and $\{z=\pi\}$ in generic position.
More precisely, none of the  values
\bel{never1}\left\{\bega{rl} (w,w_X,w_{XX}) &=~(\pi,0,0),\\[3mm]
(z,z_Y,z_{YY}) &=~(\pi,0,0),\enda\right.\eeq
\bel{never2}\left\{\bega{rl} (w,z,w_X) &=~(\pi,\pi,0),\\[3mm]
(w,z,z_Y) &=~(\pi,\pi,0),\enda\right. \eeq
\bel{never3}
\left\{\bega{rl} (w,w_X,c'(u)) &=~(\pi,0,0),\\[3mm]
(z,z_Y, c'(u)) &=~(\pi,0,0).\enda\right.\eeq
is ever attained, at any point $(X,Y)$ for which
\bel{0T}
\bigl(x(X,Y),~t(X,Y)\bigr)~\in~\R \times [0,T]\,.\eeq
}
\v
The singularities of the solution
$u$ in the $x$-$t$ plane correspond to the image of the level sets
$\{w=\pi\}$ and $\{z=\pi\}$ w.r.t.~the map  
\bel{Lambda}
\Lambda:(X,Y)~\mapsto~ \bigl(x(X,Y),~t(X,Y)\bigr).\eeq   
If none of the values
in (\ref{never1})-(\ref{never3}) is ever attained, 
by the implicit function theorem  the above level sets are the union
of a locally finite family of $\C^2$ curves in the $X$-$Y$ plane. 
In turn, restricted to the domain $\R\times [0,T]$,
the singularities of $u$ are located along finitely many
$\C^2$ curves in the $x$-$t$ plane.

\begin{figure}[htbp]
\centering
\includegraphics[width=0.7\textwidth]{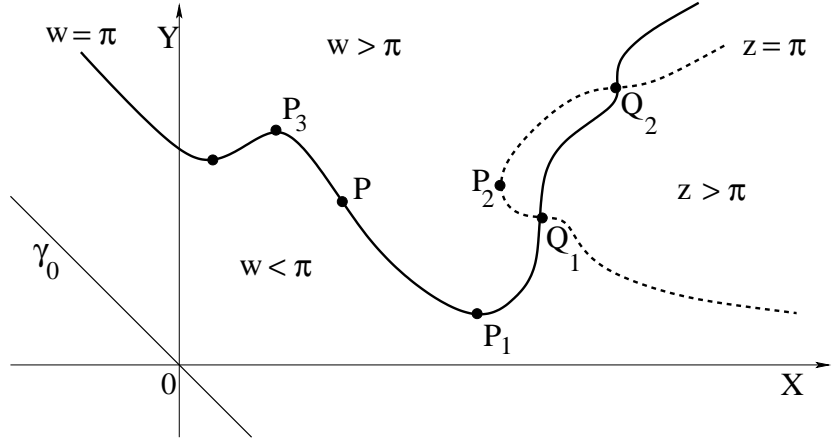}
\caption{ \small Two level sets $\{w=\pi\}$ and $\{z=\pi\}$,
in a generic conservative solution of (\ref{u})--(\ref{pq}). 
Here $P$ is a singular point of Type~1, while 
$P_1$, $P_2$, $P_3$ are  points of Type~2, and $Q_1$, $Q_2$
are  points of Type~3.  
Notice that at $P_1$, structural stability requires that 
the function $Y(X)$ implicitly defined by $w(X,Y(X))=\pi$ has strictly positive
second derivative. At the points
 $Q_1$, $Q_2$, by (\ref{wz}) one has 
 $w_Y=z_X=0$. Hence the two curves
 $\{w=\pi\}$ and $\{ z=\pi\}$ have a perpendicular intersection. }
\label{f:wa72}
\end{figure}

\begin{figure}[htbp]
\centering
\includegraphics[width=0.6\textwidth]{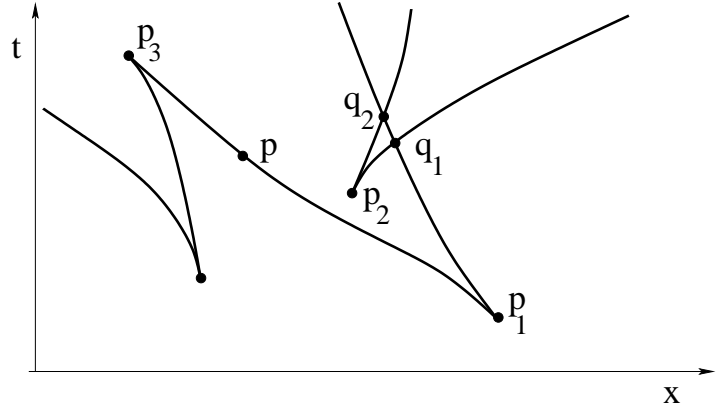}
\caption{ \small  The images of the level sets
$\{w=\pi\}$ and $\{z=\pi\}$ in Fig.~\ref{f:wa72}, under the map
$\Lambda:(X,Y)\mapsto (x(X,Y), t(X,Y))$.  In the $x$-$t$ plane, 
these represents the curves where $u=u(x,t)$ is not differentiable.
A generic solution of (\ref{1.1}) with smooth 
initial data
remains smooth outside finitely many singular points and 
finitely many singular curves, where
$u_x\to\pm\infty$.  Here $p_1,p_2, p_3$ are singular points 
where two new singular curves originate, or two singular curves 
merge and disappear.   
At the points $q_1,q_2$ a forward and a backward singular curve 
cross each other. 
}
\label{f:wa73}
\end{figure}

\section{Singularities of conservative solutions}
\setcounter{equation}{0}
For smooth  data $u_0,u_1\in\C^\infty(\R)$, the solution $(X,Y)\mapsto (x,t,u,w,z,p,q)(X,Y)$ of the semilinear system
(\ref{u})--(\ref{t}), with initial data as in (\ref{bdata}), 
(\ref{data}),  
remains smooth on the entire $X$-$Y$ plane.   
Yet, the solution
$u=u(x,t)$ of (\ref{1.1}) can have singularities because the 
coordinate change $\Lambda:(X,Y)\mapsto (x,t)$ is not smoothly invertible.
By (\ref{t})-(\ref{x}), its Jacobian matrix is computed by
\bel{J}D\Lambda~=~
\left(\bega{cc} x_X & x_Y\cr
t_X& t_Y\enda\right)~=~
\left(\bega{ccc} {(1+\cos w)\,p\over 4} && -{(1+\cos z)\,q\over 4}\\[3mm]
{(1+\cos w)\,p\over 4c(u)}&& {(1+\cos z)\,q\over 4c(u)}\enda\right)
\eeq
We recall that $p,q$ remain uniformly positive and uniformly bounded
on compact subsets of the $X$-$Y$ plane.
By Remark~1, at a point $(X_0,Y_0)$ where  
$w\not=\pi$ and $z\not=\pi$, this  matrix
is  invertible, having a strictly positive determinant.   
The function $u=u(x,t)$ considered at 
(\ref{graph}) is thus smooth on a neighborhood of the point 
$$(x_0, t_0) ~=~\bigl( x(X_0, Y_0)\,,~t(X_0,Y_0)\Bigr).$$
To study the set of points $x$-$t$ plane where $u$ is singular,
we thus need to look at points where either $w=\pi$ or $z=\pi$.

If the generic conditions (\ref{never1})--(\ref{never3}) are satisfied,
then we have the implications
$$\left\{\bega{rl}w~=~\pi\quad\hbox{and}\quad  w_X~=~0\qquad &\ds \implies\qquad w_Y~=~
{c'(u)\over 8c^2(u)}(\cos z+1)q~\not=~0\,,\\[4mm]
z~=~\pi\quad\hbox{and}\quad  z_Y~=~0\qquad &\ds \implies\qquad z_X~=~
{c'(u)\over 8c^2(u)}(\cos w+1)p~\not=~0\,.\enda\right.$$
Therefore, by the implicit function theorem,
the level sets 
\bel{Swz}\S^w~\doteq~\{(X,Y)\,;~~w(X,Y)=\pi\}\,,\qquad\qquad
\S^z\doteq\{(X,Y)\,;~~z(X,Y)=\pi\}\,,\eeq
are the union of a locally finite family of smooth curves.
The singularities of $u$ 
in the $x$-$t$ plane are contained in the images of $S^w$ and $S^z$
under the map (\ref{Lambda}).
Relying on Theorem~2, 
we shall distinguish three types of singular points $P=(X_0,Y_0)$.
\v
\begi
\item[{\bf (1)}]  Points  where $w=\pi$ but $w_X\not= 0$ and $z\not= \pi$
(or else, where   $z=\pi$ but $z_Y\not= 0$ and $w\not= \pi$).

\item[{\bf (2)}] Points  where $w=\pi$ and  $w_X= 0$, 
but  $w_{XX}\not= 0$
(or else: $z=\pi$ and  $z_Y= 0$, but  $z_{YY}\not= 0$).   

\item[{\bf (3)}]  Points  where $w=\pi$ and $z=\pi$. 
\endi

Points of Type 1 form a locally finite family of 
$\C^2$ curves in the $X$-$Y$ plane (Fig.~\ref{f:wa72}).   
Their images $\Lambda(P)$
yield a family of characteristic curves in the $x$-$t$ plane where the solution $u=u(x,t)$ is 
singular (i.e., not differentiable).  

Points of Type 2 are isolated. Their images in the $x$-$t$ plane are points
where two singular curves 
 initiate or terminate (Fig.~\ref{f:wa73}).

Points of Type 3 are those where two curves $\{w=\pi\}$ and $\{z=\pi\}$ intersect.
Their image in the $x$-$t$ plane are points where two singular curves
cross, 
with speeds $\pm c(u)$.
\v
Our main result provides a detailed description of the 
solution $u=u(x,t)$
in a neighborhood of each one of these singular points.  For simplicity, we shall 
assume that the initial data $(u_0,u_1)$ in (\ref{1.2}) are smooth,
so we shall not need to count how many derivatives are actually used to 
derive the Taylor approximations.
\v
{\bf Theorem 3.} {\it  Let the assumptions {\bf (A)} hold, and consider 
generic initial data  $(u_0,u_1)\in\D$ as in (\ref{U}), with $u_0, u_1\in C^\infty(\R)$.   
Call  $(u,w,z,p,q,x,t)$  the corresponding solution  of the semilinear 
system (\ref{u})--(\ref{t}) and let 
$u=u(x,t)$ be the  solution to the original equation (\ref{1.1}).
Consider a singular point $P=(X_0,Y_0)$ where $w=\pi$, 
and set $(x_0, t_0)\doteq (x(X_0,Y_0), t(X_0,Y_0))$.
\begi
\item[(i)] If $P$ is a point of Type 1, along a curve where $w=\pi$, 
then there exist constants
$a\not= 0$ and $b_1,b_2$ such that
\bel{T1}\bega{rl}
u(x,t)&=~u(x_0, t_0) - a\cdot  \Big[c(u_0)(t-t_0)+(x-x_0)\Big]^{2/3} \cr\cr
&\ds+ b_1\cdot(x-x_0) + b_2\cdot (t-t_0)+
\O(1)\cdot \Big(|t-t_0|+ |x-x_0|\Big)^{4/3 }.
\enda
\eeq

\item[(ii)] If $P$ is a point of Type 2, where $w=\pi$,  $w_X=0$, and
$w_{XX}>0$,
then in the $x$-$t$ plane this corresponds to a point $(x_0, t_0)$
where two new singular curves $\gamma^-, \gamma^+$ originate.
In this case, there exists a constant  $a\not= 0$  such that
\bel{T2}\bega{rl}
u(x,t)&=~u(x_0, t_0) + a\cdot  \Big[c(u_0)(t-t_0)+(x-x_0)\Big]^{3/5}
 + 
\O(1)\cdot \Big(|t-t_0|+ |x-x_0|\Big)^{4/5}.\enda\eeq

\item[(iii)] 
 If $P$ is a point of Type 3, where $w=z=\pi$,
then in the $x$-$t$ plane this corresponds to a point $(x_0, t_0)$
where two  singular curves $\gamma, \tilde \gamma$ cross each other.
In this case, there exist constants  $a_1\not= 0$ and 
$a_2\not= 0$ such that
\bel{T3}\bega{rl}
u(x,t)&=~u(x_0, t_0) + a_1\cdot  
\Big[c(u_0)(t-t_0)+(x-x_0)\Big]^{2/3}\cr\cr
&\qquad+ a_2\cdot  
\Big[c(u_0)(t-t_0)-(x-x_0)\Big]^{2/3} + 
\O(1)\cdot \Big(|t-t_0|+ |x-x_0|\Big).\enda\eeq
\endi
}
\v
Throughout  the following, given
 a point $P=(X_0, Y_0)$ in the $X$-$Y$ plane
where $w=\pi$,  we 
denote by $(u_0, w_0, z_0, p_0, q_0,x_0, t_0)$
the values of $(u,w,z,p,q,x,t)$ at $(X_0, Y_0)$.
The three parts of Theorem~3 will be proved separately.

\subsection{Singular curves.}
Let $P=(X_0, Y_0)$ be a point of Type 1,  where
\bel{Case1}
w_0~=~\pi,\qquad\qquad z_0~\not=~ \pi,\qquad
\qquad w_X(X_0, Y_0)~\not=~ 0.\eeq
By the implicit function theorem, the level set where 
$w=\pi$ is locally the graph of a smooth function
$X=\Phi(Y)$, with $\Phi(Y_0)=X_0$.
We claim that, in a neighborhood of the point 
$(x_0, t_0)= \Lambda(X_0, Y_0)$, the image $\Lambda(\S^w)$ 
 is a smooth curve in the $x$-$t$ plane, say 
 \bel{gamma}\gamma~=~\bigl\{(x,t)\,;~
 x=\phi(t)\bigr\}.\eeq
Indeed, the curve $\gamma$ is the image of the smooth 
curve $\{X=\Phi(Y)\}$
under the smooth, one-to-one map
$$Y~\mapsto~\bigl(x(\Phi(Y),Y)\,,~ t(\Phi(Y),Y)\bigr).$$
For future record, 
we compute the first two derivatives of $\phi$ at $t=t_0$.
Differentiating   the identity
$w(\Phi(Y),Y)=\pi$
one obtains
$$w_X \Phi' + w_Y ~=~0\,,$$
$$w_{XX}\cdot (\Phi')^2 + 2 w_{XY} \Phi' + w_{YY} + w_X \Phi''~=~0\,.$$
By (\ref{x})-(\ref{t}),  at the point $(X_0,Y_0)$ we have
$${d\over dY}\bigl(x(\Phi(Y),Y)\,,~ t(\Phi(Y),Y)\bigr)~=~
\left( -{(1+\cos z_0)\,q_0\over 4}\,, ~{(1+\cos z_0)\,q_0\over 4c(u_0)}
\right)~\not=~(0,0).
$$
Observing that
$$\phi'(t(\Phi(Y),Y))~=~{x_X(\Phi(Y),Y)\cdot\Phi' (Y)+ x_Y(\Phi(Y),Y)
\over t_X(\Phi(Y),Y)\cdot \Phi'(Y) + t_Y(\Phi(Y),Y)}\,,$$
at $t=t_0$ we have
$$\phi'(t_0)~=~ -c(u_0)\,.$$
In a similar way we find
$$\phi''(t_0)~=~\frac{x_{YY}(X_0, Y_0)t_Y(X_0, Y_0) - t_{YY}(X_0, Y_0)x_Y(X_0, Y_0)}{t_Y^3(X_0, Y_0)}~=~-\frac{c'(u_0)\sin z_0}{1+\cos z_0}\,.$$

\v
Next, by (\ref{u}) one has  
\bel{uY} u_X(X_0, Y_0)~=~0\,,\qquad\qquad 
u_Y(X_0, Y_0) ~=~ \frac{\sin z_0}{4c(u_0)} q_0~\doteq~\alpha_1\,.\eeq
Differentiating the first equation in (\ref{u})
w.r.t.~$X$ and using (\ref{wz})-(\ref{pq}) we obtain
$$u_{XX}~=~{\cos w\over 4c(u)} w_X p - 
{\sin w\over 4 c^2(u)} c'(u) \cdot
{\sin w\over 4c(u)} p^2 + {\sin w\over 4c(u)} p_X\,,$$
\bel{uXX}
u_{XX}(X_0, Y_0)~=~{w_X(X_0, Y_0)\over 4c(u_0)}p_0~\doteq~\alpha_2
~\not=~0
\,,\eeq
\bel{uXXX} u_{XXX}(X_0, Y_0)~=~-{1\over 4c(u_0)} 
\Big(w_{XX}(X_0, Y_0)~ p_0 +
2 w_X (X_0, Y_0)p_X(X_0, Y_0)\Big)~\doteq~\alpha_3\,,\eeq
\bel{uXY} u_{XY}(X_0, Y_0)~=~-{p_0\over 4 c(u_0)} \cdot {c'(u_0)
\over 8c^2(u_0)} (\cos z_0+1) q_0~\doteq~\alpha_4\,.\eeq
This yields the local Taylor approximation
\bel{up1}
\bega{rl}
u(X, Y)&= ~\ds u_0 +  \alpha_1\, (Y-Y_0) 
+{\alpha_2\over 2} \,(X-X_0)^2 + {\alpha_3\over 6}\, (X-X_0)^3
+\alpha_4\,(X-X_0)(Y-Y_0)\cr\cr
&\qquad\ds + \O(1)\cdot \Big(|X-X_0|^4+|Y-Y_0|^2+
|X-X_0|^2\,|Y-Y_0|\Big).
\enda
\eeq
Using (\ref{t}), we perform an entirely similar computation 
 for the function $t$ in a neighborhood
of $(X_0, Y_0)$.
\bel{t1}t_X(X_0, Y_0) ~=~0,\qquad\qquad
t_Y(X_0, Y_0)
~=~{1+\cos z_0\over 4c(u_0)}q_0~\doteq~\beta_1~>~ 0\,,\eeq
$$
t_{XX}~=~-{\sin w\over 4c(u)} w_X p-{1+\cos w\over 4c^2(u)} c'(u)\,
u_X p + {1+\cos w\over 4c(u)}\, p_X\,,$$
\bel{t222}t_{XX}(X_0, Y_0)~=~t_{XY}(X_0, Y_0)~=~0\,,\eeq
\bel{t3}
t_{XXX}(X_0, Y_0)~=~{w_X^2(X_0, Y_0)\over 4c(u_0)}\, p_0~\doteq~\beta_3~\not=~0\,.
\eeq
This yields the Taylor approximation
\bel{tXY}\bega{rl}
t(X, Y)& = ~\ds 
t_0 +\beta_1\,(Y-Y_0) + {\beta_3\over 6}\,(X-X_0)^3\cr\cr
&\quad+ \O(1)\cdot \Big(|X-X_0|^4+|Y-Y_0|^2+|X-X_0|^2\,|Y-Y_0|\Big).
\enda\eeq
Finally, for the function $x$, using (\ref{x}) we find
\bel{x1}
 x_X(X_0, Y_0)~=~0\,,\qquad\qquad 
 x_Y(X_0, Y_0)~=~-{1+\cos z_0\over 4}\,q_0~\doteq~-\gamma_1~<~0\,,
 \eeq
$$
x_{XX}~=~-{\sin w\cdot w_X\over 4}  p + {1+\cos w\over 4}\, p_X\,,$$

\bel{x2}
x_{XX}(X_0, Y_0)~=~0,\qquad\qquad x_{XY}(X_0, Y_0)~=~0,\eeq
\bel{x3}
x_{XXX}(X_0, Y_0)~=~{w_X^2(X_0, Y_0)\over 4}\, p_0~\doteq~\gamma_3~>~0.
\eeq
This yields the Taylor approximation
\bel{xXY}\bega{rl}x(X, Y)& =\ds ~x_0 -\gamma_1\,(Y-Y_0) 
+ \gamma_3\,(X-X_0)^3 \cr\cr
&\quad+  \O(1)\cdot \Big(|X-X_0|^4+|Y-Y_0|^2+|X-X_0|^2\,|Y-Y_0|\Big).
\enda
\eeq
Observing that the above Taylor coefficients satisfy
\bel{bg}
\gamma_1~=~c(u_0)\,\beta_1\,,\qquad\qquad \gamma_3~=~c(u_0)\,\beta_3\,,
\eeq
from (\ref{tXY}) and (\ref{xXY}) we deduce 
\bel{99}\bega{l}
(x-x_0)- c(u_0)(t-t_0)\\[3mm]
\qquad =~-2\gamma_1\,(Y-Y_0)+
\O(1)\cdot \Big(|X-X_0|^4+|Y-Y_0|^2+|X-X_0|^2\,|Y-Y_0|\Big),\\[4mm]
(x-x_0)+ c(u_0)(t-t_0)\\[3mm]
\qquad =~2\gamma_3\,(X-X_0)^3+
\O(1)\cdot \Big(|X-X_0|^4+|Y-Y_0|^2+|X-X_0|^2\,|Y-Y_0|\Big)\,.\enda
\eeq

Next,
using (\ref{tXY}) and (\ref{xXY}) we obtain  an approximation 
for $X,Y$ in terms of $x,t$, namely
$${1+ \cos z_0\over 2} q_0 (Y-Y_0)~=~c(u_0)(t-t_0)- (x-x_0) 
+ \O(1)\cdot \Big(|X-X_0|^4+|Y-Y_0|^2+|X-X_0|^2\,|Y-Y_0|\Big).$$
$${w_X^2(X_0, Y_0)\over 12} p_0 (X-X_0)^3~=~c(u_0)(t-t_0)+ (x-x_0) 
+ \O(1)\cdot \Big(|X-X_0|^4+|Y-Y_0|^2+|X-X_0|^2\,|Y-Y_0|\Big).$$
Inserting the two above expressions into \eqref{up1}, we finally obtain
\bel{utx1}\bega{rl}
u(t, x)& = \ds~u(t_0, x_0) - \left( {9p_0\over 32 \, w_X(X_0, Y_0)}\right)^{1/3} \Big[c(u_0)(t-t_0)+(x-x_0)\Big]^{2/3} \cr\cr
&\qquad \ds+ \,\frac{\sin z_0}{2c(u_0)(1+\cos z_0)}\Big[c(u_0)(t-t_0) - (x-x_0)\Big] \cr\cr
&\qquad \ds-\frac{w_{XX}(X_0, Y_0)}{2c(u_0)w_X^2(X_0, Y_0)}\Big[c(u_0)(t-t_0)+(x-x_0)\Big]\cr\cr
 &\qquad \ds+\,\O(1)\cdot \Big(|t-t_0|+ |x-x_0|\Big)^{4/3 }.
\enda
\eeq
This proves (\ref{T1}), with 
\bel{ab12}
a~=~ \left( {9p_0\over 32 \, w_X(X_0, Y_0)}\right)^{1/3}~\not=~0\,.
\eeq
The coefficients $b_1,b_2$ can also 
be easily computed from (\ref{utx1}).

\begin{figure}[htbp]
\centering
\includegraphics[width=0.9\textwidth]{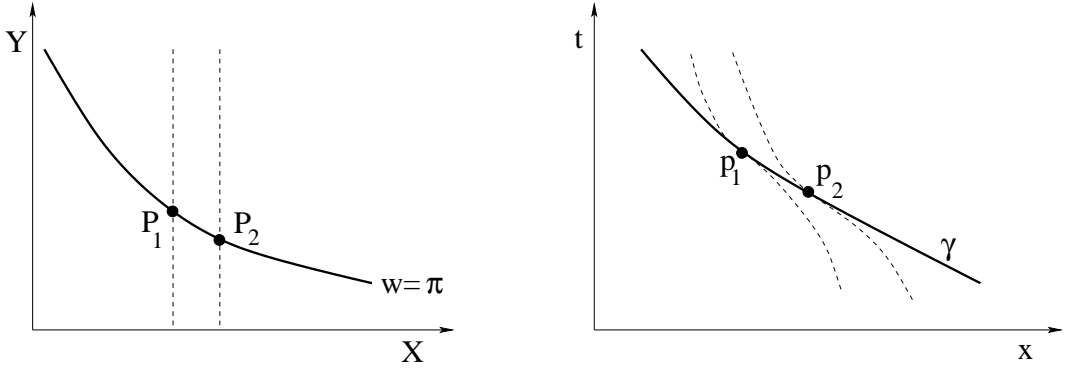}
\caption{ \small  Left: a singular curve where $w=\pi$, in the $X$-$Y$ plane. 
Vertical lines where $X$ = constant correspond to characteristic curves
of the wave equation (\ref{1.1}), where $\dot x = -c(u)$.
Right: the images of these curves in the $x$-$t$ plane, under 
the map $\Lambda$ at (\ref{Lambda}). 
The singular curve $\gamma$ is an envelope of characteristic curves, 
which cross it tangentially. }
\label{f:wa74}
\end{figure}

{\bf Remark 3.} By (\ref{T1}), the solution $u$ is only H\"older continuous
of exponent 2/3 near the singular curve $\gamma$ in (\ref{gamma}).
In particular, the Cauchy problem
$$\dot x(t)~=~-c\bigl(u(t,x(t))\bigr),\qquad\qquad
x(t_0)~=~ \phi(t_0)\,,$$
has a solution $t\mapsto x(t)$ which crosses $\gamma$ 
at the point $(x_0, t_0)$. 
Calling $\delta(t)\doteq x(t)-\phi(t)$,
to leading order one has
$$\dot \delta~=~ c'(u_0)\cdot a\,\delta^{2/3}\,.$$
Hence, for $t\approx t_0$ we have
\bel{char}\delta(t) ~\approx~\left({c'(u_0)\, a\over 3}\right)^3
(t-t_0)^3.\eeq
 The singular curve $\gamma$ is thus an envelope of characteristic curves, which cross it tangentially (see Fig.~\ref{f:wa74}).
 
\v
\subsection{Points where two  singular curves originate or terminate.}
Let $P=(X_0, Y_0)$ be a point of Type~2,
where
\bel{sp2}w_0~=~\pi,\qquad z_0~\not= ~\pi,\qquad 
 w_X(X_0, Y_0)~=~0, \qquad w_{XX}(X_0, Y_0)~\not=~0\,.
\eeq

Recalling (\ref{impl}),
by \eqref{sp2} we have
$$
w_Y(X_0, Y_0) ~= ~\frac{c'(u_0)}{8c^2(u_0)}(1+\cos z_0) q_0~\neq ~0.$$
By (\ref{u}), at the point $(X_0, Y_0)$ we have
$$
u_X ~= ~u_{XX} ~=~ 0\,,\qquad\qquad \, u_Y~ = ~\frac{\sin z_0}{4c(u_0)}\,q_0\,,$$
$$
u_{XY} ~= ~-\frac{c'(u_0)}{32c^3(u_0)}(1+\cos z_0) p_0 q_0\,,
\qquad\qquad  u_{XXX}~ =~ -\frac{w_{XX}(X_0, Y_0)}{4c(u_0)} p_0\,.
$$
In this case, the Taylor approximation for $u$ near the point $(X_0,Y_0)$ takes the form
\bel{up2}
\bega{rl}
u(X, Y) &=~\ds u_0 + \frac{\sin z_0}{4c(u_0)}q_0\,(Y-Y_0) - \frac{w_{XX}(X_0, Y_0)}{24c(u_0)}p_0\,(X-X_0)^3\cr\cr
& \qquad\ds
+ \O(1)\cdot\Big(|X-X_0|^4 +|Y-Y_0|^2 + |X-X_0|\,|Y-Y_0|\Big).
\enda
\eeq 
Computing the partial derivatives of $x(X,Y)$ and $t(X,Y)$ 
at the point $(X_0, Y_0)$,
by (\ref{wz}) and (\ref{sp2}) we find 
\bel{xX}x_X~=~x_{XX}~=~x_{XXX}~=~x_{XXXX}~=~x_{XY}~=~x_{XXY}~=~0,\eeq
\bel{xXX}x_{XXXXX}~=~{3w_{XX}^2(X_0,Y_0)\over 4}\, p_0~\not=~0\,,\qquad\qquad 
x_Y~=~-{1+\cos z_0\over 4}\, q_0~\not=~0\,.\eeq
\bel{tX}t_X~=~t_{XX}~=~t_{XXX}~=~t_{XXXX}~=~t_{XY}~=~t_{XXY}~=~0,\eeq
\bel{tXX}t_{XXXXX}~=~{3w_{XX}^2(X_0,Y_0)\over 4c(u_0)}\, p_0~\not=~0\,,\qquad\qquad 
t_Y~=~{1+\cos z_0\over 4c(u_0)}\, q_0~\not=~0\,.\eeq
This yields the Taylor approximations
\bel{x22}
\bega{rl}
x(X, Y) &=~ \ds x_0 - \frac{1+\cos z_0}{4}q_0(Y-Y_0) + 
\frac{3w_{XX}^2(X_0, Y_0)}{5!\,4}p_0(X-X_0)^5\\[4mm]
&\qquad\ds + \O(1)\cdot \Big(|X-X_0|^6 + |Y-Y_0|^2\Big).
\enda
\eeq
\bel{t22}
\bega{rl}
t(X, Y) &=~ \ds t_0 +
\frac{1+\cos z_0}{4c(u_0)}q_0(Y-Y_0) +
\frac{3w_{XX}^2(X_0, Y_0)}{5!\,4 \,c(u_0)}p_0(X-X_0)^5\\[4mm]
&\qquad\ds+ \O(1)\cdot \Big(|X-X_0|^6 + |Y-Y_0|^2 \Big),\enda
\eeq
Combining (\ref{x22}) with (\ref{t22}) we obtain
\bel{e33}\bega{rl}(X-X_0)^5 &=~\ds
\frac{5!\,2 }{3w_{XX}^2(X_0, Y_0)\, p_0}
\cdot \bigl[c(u_0)(t-t_0) +(x-x_0)\bigr]
\cr\cr&\qquad\qquad
+ \O(1)\cdot \Big(|X-X_0|^6 + |Y-Y_0|^2 \Big),\enda\eeq
\bel{e44}\bega{rl}Y-Y_0 &=~\ds
\frac{2}{(1+\cos z_0)\,q_0}\cdot 
\bigl[c(u_0)(t-t_0) -(x-x_0)\bigr]
+ \O(1)\cdot \Big(|X-X_0|^6 + |Y-Y_0|^2 
\Big).\enda\eeq
Inserting (\ref{e33})-(\ref{e44}) into \eqref{up2} 
we eventually obtain
\bel{utx2}\bega{rl}
u(t, x)& 
=\ds~ u(t_0, x_0) - \frac{1}{24c(u_0)}\cdot \left(\frac{80^3\, p_0^2}{
w_{XX}(X_0, Y_0)}\right)^{1/5}\cdot
\Big[c(u_0)(t-t_0) + (x-x_0)\Big]^{3/5} \cr\cr
&
\qquad+ \O(1)\cdot \Big(|t-t_0|+ |x-x_0|\Big)^{4/5}.
\enda
\eeq
This proves (\ref{T2}).

It remains to show that two 
 singular curves originate or terminate at the point $(x_0, t_0)$.
 To fix the ideas, assume that
\bel{kdef}\kappa~\doteq~-{w_{XX}(X_0, Y_0)\over 2 w_Y(X_0, Y_0)}~>~0
\,.\eeq
By the implicit function theorem,
the curve where $w=\pi$ can be approximated as
\bel{wpi}
Y-Y_0~=~\kappa (X-X_0)^2 + \O(1)\cdot |X-X_0|^3.\eeq
On the other hand, by  (\ref{t22}) we have
\bel{Yt}Y-Y_0~=~\alpha(t-t_0) - \beta(X-X_0)^5 + \O(1)\cdot \Big(|X-X_0|^6 + |t-t_0|^2 + |X-X_0|\,|t-t_0|\Big),
\eeq
with 
$$\alpha~=~{4c(u_0)\over (1+\cos z_0)q_0}~>~0\,,
\qquad\qquad 
\beta~=~{3w^2_{XX}(X_0, Y_0)\, p_0\over 5! \,(1+\cos z_0)q_0}~>~0\,.$$
Combining  (\ref{wpi}) with  (\ref{Yt}) we obtain
\bel{roots}
\kappa(X-X_0)^2~=~\alpha(t-t_0) + \O(1)\cdot |X-X_0|^3.\eeq
Therefore, as shown in Fig.~\ref{f:wa75} in a neighborhood
of $(X_0, Y_0)$ the following holds:
\begi
\item The two curves $\{t(X,Y)=t_0\}$ and $\{w(X,Y)=\pi\}$
intersect exactly at the point $(X_0, Y_0)$.
\item When $\tau<t_0$, the curves $\{t(X,Y)=\tau\}$ and $\{w(X,Y)=\pi\}$
have no intersection.
\item When $\tau>t_0$, the curves $\{t(X,Y)=\tau\}$ and $\{w(X,Y)=\pi\}$
have two intersections, at points 
$P_1=(X_1,Y_1)$ and $P_2=(X_2, Y_2)$ with 
\bel{es5}
\left\{\bega{rl}X_1-X_0&=\ds ~-\sqrt{{\alpha\over\kappa}(\tau-t_0)} + \O(1)\cdot (\tau-t_0),\\[4mm]
X_2-X_0&=\ds ~+\sqrt{{\alpha\over\kappa}(\tau-t_0)} + \O(1)\cdot  (\tau-t_0).\enda\right.
\eeq
\bel{es6}\left\{\bega{rl}Y_1-Y_0&= ~\alpha(\tau-t_0) + \O(1)\cdot (\tau-t_0)^{3/2},\cr
Y_2-Y_0&= ~\alpha(\tau-t_0) + O(1)\cdot (\tau-t_0)^{3/2},\enda\right.\eeq
\endi

\begin{figure}[htbp]
\centering
\includegraphics[width=1.0\textwidth]{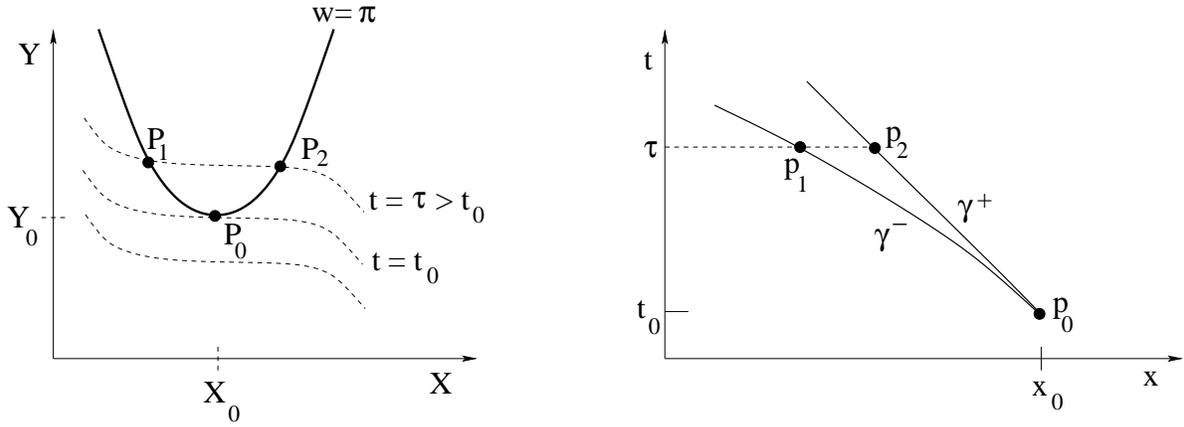}
\caption{ \small  Left: the equation $w(X,Y)=\pi$ implicitly 
defines a function $Y(X)$ with a strict local minimum at $X_0$.  Under generic conditions, $Y''(X_0)>0$.  
The dotted curves where $t(X,Y) = \tau$ have  0, 1, or 2 intersections
respectively,
if $\tau<t_0$, $\tau=t_0$, or $\tau>t_0$. 
Right: the image of the curve $\{w=\pi\}$ under the map $\Lambda$ in (\ref{Lambda})
consists of two singular curves $\gamma^-, \gamma^+$ starting at the point 
$p_0=(x_0, t_0)$.   For $\tau>t_0$, the distance between these two curves is 
$\gamma^+(\tau)-\gamma^-(\tau)=\O(1)\cdot (\tau-t_0)^{5/2}$.
 }
\label{f:wa75}
\end{figure}
For $t>t_0$, the solution $u=u(t,x)$ is thus singular along two
curves $\gamma^-, \gamma^+$ in the $x$-$t$ plane
(see Fig.~\ref{f:wa75}, right). Our next goal
is to derive an asymptotic description of these curves in
a neighborhood of the point $(x_0, t_0)$, namely
\bel{g12}
\left\{ \bega{rl}
\gamma^-(t)&=~x_0 - c(u_0)(t-t_0) +\tilde \alpha (t-t_0)^2 - \tilde\beta (t-t_0)^{5/2} + 
 \O(1)\cdot(t-t_0)^3,\\[3mm]
\gamma^+(t)&=~x_0 - c(u_0)(t-t_0) + \tilde \alpha (t-t_0)^2 + \tilde \beta (t-t_0)^{5/2} + 
 \O(1)\cdot(t-t_0)^3,\enda
\right.\eeq
for suitable constants $\tilde\alpha,\tilde\beta$.

To prove (\ref{g12}), we need to compute
 more accurate Taylor approximations for  $t$ and $x$
 near the point $(X_0, Y_0)$. 
\bel{tt2}
\bega{rl}
t(X, Y) &\ds=~ t_0 + \frac{1+\cos z_0}{4c(u_0)}q_0(Y-Y_0)  +a(Y-Y_0)^2 + 
\frac{3w_{XX}^2(X_0, Y_0)}{5!\,4c(u_0)}p_0(X-X_0)^5\cr\cr
&\qquad \ds+\frac{w_{Y}^2(X_0, Y_0)}{8c(u_0)}p_0(X-X_0)(Y-Y_0)^2\cr\cr
&\qquad
+ \O(1)\cdot \Big(|X-X_0|^6 + |Y-Y_0|^3 + |X-X_0|^2\,|Y-Y_0|^2\Big),
\enda
\eeq
 \bel{xx2}
\bega{rl}
x(X, Y)& \ds=~ x_0 - \frac{1+\cos z_0}{4}q_0(Y-Y_0) +b(Y-Y_0)^2+ 
\frac{3w_{XX}^2(X_0, Y_0)}{5!\,4}p_0(X-X_0)^5\cr\cr
&\ds\qquad +\frac{w_{Y}^2(X_0, Y_0)}{8}p_0(X-X_0)(Y-Y_0)^2 
\cr\cr
&\qquad \ds + \O(1)\cdot \Big(|X-X_0|^6 + |Y-Y_0|^3 + |X-X_0|^2\,|Y-Y_0|^2\Big).
\enda
\eeq
The constants $a$ and $b$ are here given by
$$a~ = ~-\frac{z_Y\sin z}{8c(u)}q  + \frac{1+\cos z}{8c(u)}\left(q_Y - \frac{c'(u)\sin z}{4c^2(u)}q^2\right),$$
$$b~ =~ \frac 18z_Yq\sin z  - \frac 18(1+\cos z)q_Y\,.$$
where the right hand sides are evaluated at the point $(X_0, Y_0)$.

For a fixed $\tau>t_0$, let $P_1=(X_1,Y_1)$ and $P_2=(X_2,Y_2)$ be 
 the two points  where the curves $\{t(X,Y)=\tau\}$ and 
 $\{w(X,Y)=\pi\}$ intersect.  
 Let $x=\gamma^-(\tau)$ and $x=\gamma^+(\tau)$
 describe the 
 corresponding points in the $x$-$t$ plane (see Fig.~\ref{f:wa75}). 

At the  intersection point $P_1=(X_1,Y_1)$, 
using \eqref{es5} and \eqref{es6} we obtain
\bel{gpf12}
\bega{l}
x(X_1, Y_1)-x_0+c(u_0)(\tau-t_0) \ds\cr\cr
 =~\ds (a+b)(Y_1-Y_0)^2+ 
\frac{3w_{XX}^2(X_0, Y_0)}{5!\,2}p_0(X_1-X_0)^5
+\frac{w_{Y}^2(X_0, Y_0)}{4}p_0(X_1-X_0)(Y_1-Y_0)^2\cr\cr
\qquad \ds + \O(1)\cdot \Big(|X_1-X_0|^6 + |Y_1-Y_0|^3 + 
|X_1-X_0|^2\,|Y_1-Y_0|^2\Big)\cr\cr

 =\ds ~\alpha^2(a+b)(\tau-t_0)^2-\left( 
\frac{3w_{XX}^2(X_0, Y_0)}{5!\,2\,\kappa^{5/2}}+\frac{w_{Y}^2(X_0, Y_0)}
{4\,\kappa^{1/2}}\right)
\alpha^{5/2}p_0(\tau-t_0)^{5/2}+ \O(1)\cdot (\tau-t_0)^3.
\enda
\eeq
This yields
the equation for  $\gamma^-$ in \eqref{g12}, with suitable coefficients $\tilde\alpha,\tilde\beta$. 
An entirely similar argument yields the
equation for  $\gamma^+$.
In particular, the distance between these two singular curves is
\bel{dx12}
\gamma^+(t)- \gamma^-(t)~=~ 2\tilde 
\beta(t-t_0)^{5/2} + \O(1)\cdot|t-t_0|^3.
\eeq
\v
\subsection{Points where two singular curves cross.}
We now consider a point $P=(X_0, Y_0)$ where
$ w=z = \pi$.

For a generic solution, satisfying the conclusion of Theorem~2, this
implies
\bel{wzx}
w_X(X_0, Y_0) ~\neq~ 0, \qquad \qquad 
z_Y(X_0, Y_0)~ \neq~ 0\,.
\eeq
On the other hand,  \eqref{wz} yields
$$w_Y (X_0, Y_0)~=~ z_X (X_0, Y_0)~=~0.$$
By (\ref{u}) and (\ref{cc1}), we know that
$$ u_X(X_0, Y_0)~ =~ u_Y(X_0, Y_0) ~=~ u_{XY}(X_0, Y_0) ~=~0.$$
Hence, in a neighborhood of   $(X_0,Y_0)$ the function $u$ 
can be approximated by
\bel{up3}
u(X, Y)~ = ~ u_0 - \frac{w_X(X_0, Y_0)}{8c(u_0)}p_0(X-X_0)^2 - \frac{z_Y(X_0, Y_0)}{8c(u_0)}q_0(Y-Y_0)^2 + \O(1)\cdot \Big(|X-X_0|+ |Y-Y_0|\Big)^3.
\eeq


In addition, by (\ref{x})-(\ref{t}) we have
\bel{t33}
\bega{rl}
t(X, Y) &\ds=~ t_0 + \frac{w_X^2(X_0, Y_0)}{24c(u_0)} p_0(X-X_0)^3 + \frac{z_Y^2(X_0, Y_0)}{24c(u_0)}q_0(Y-Y_0)^3\cr\cr
&\qquad + \O(1)\cdot \Big( |X-X_0|+ |Y-Y_0|\Big)^4
,\enda\eeq
\bel{x33}\bega{rl}
x(X, Y) &\ds= x_0 + \frac{w_X^2(X_0, Y_0)}{24} p_0(X-X_0)^3 - \frac{z_Y^2(X_0, Y_0)}{24}q_0(Y-Y_0)^3\cr\cr
&\qquad + \O(1)\cdot \Big( |X-X_0|+ |Y-Y_0|\Big)^4.
\enda\eeq

Using (\ref{t33})-(\ref{x33}) in (\ref{up3}) we eventually obtain
\bel{utx3}
\bega{rl}
u(t, x) &=\ds 
u(t_0, x_0) -\frac{1}{8c(u_0)}
\left(\frac{144 \,p_0}{w_X(X_0, Y_0)}\right)^{1/3} \bigl[c(u_0)(t-t_0) + (x-x_0)\bigr]^{2/3}\cr\cr
&\ds -\frac{1}{8c(u_0)}
\left(\frac{144 \,q_0}{z_Y(X_0, Y_0)}\right)^{1/3} \bigl[c(u_0)(t-t_0) - (x-x_0)\bigr]^{2/3}+ \O(1)\cdot \Big( |t-t_0| + |x-x_0|\Big).
\enda
\eeq
This proves (\ref{T3}).\endproof

\section{Dissipative solutions}
\setcounter{equation}{0}

In this last section we assume $c'(u)>0$ and 
study the structure of a dissipative solution in a neighborhood
of a point where a new singularity appears.
We recall that dissipative solutions can be characterized 
by the property that $R,S$ in (\ref{2.1}) are bounded below,
on any compact subset of the domain $\{(t,x)\,;~~t>0,~x\in\R\}$.
As proved in \cite{BH}, dissipative solutions can  
be constructed by the
same transformation of variables as in (\ref{wzdef}), (\ref{XY}), and (\ref{pqdef}).
However, the equations (\ref{wz})-(\ref{pq}) should now be replaced by
\bel{4.1}
\left\{
\bega{rl}
w_Y&=~\theta\cdot{c'(u)\over 8c^2(u)}\,( \cos z - \cos w)\,q\,,\cr\cr
z_X&=\theta\cdot{c'(u)\over 8c^2(u)}\,( \cos w - \cos z)\,p\,,\enda\right.\eeq
\bel{4.2}
\left\{\bega{rl}p_Y&=\theta\cdot{c'(u)\over 8c^2(u)}\,
\big[\sin z-\sin w\big]\,pq\,,\cr\cr
q_X&=\theta\cdot{c'(u)\over 8c^2(u)}\,
\big[\sin w-\sin z\big]\,pq\,,\enda\right.\eeq
where
\bel{thdef}
\theta~=~\left\{\bega{rl} 1\quad\hbox{if}\quad \max\{w,z\}<\pi\,,\cr
0\quad\hbox{if}\quad \max\{w,z\}\geq\pi\,.\enda\right.\eeq
Notice that, by setting  $\theta\equiv 1$, one would again recover
the conservative solutions.

It is interesting to compare a conservative and a dissipative solution, 
with the same initial data.
Consider a point $P=(X_0, Y_0)$ of Type~2, where two new singular 
curves $\gamma^-,\gamma^+$ originate, in the conservative solution.
To fix the ideas, assume that the singularity occurs in backward moving waves, so that $R\to +\infty$ but $S$ remains bounded. 
Moreover, let the conditions 
(\ref{sp2}) and (\ref{kdef}) hold.

\begin{figure}[htbp]
\centering
\includegraphics[width=1.0\textwidth]{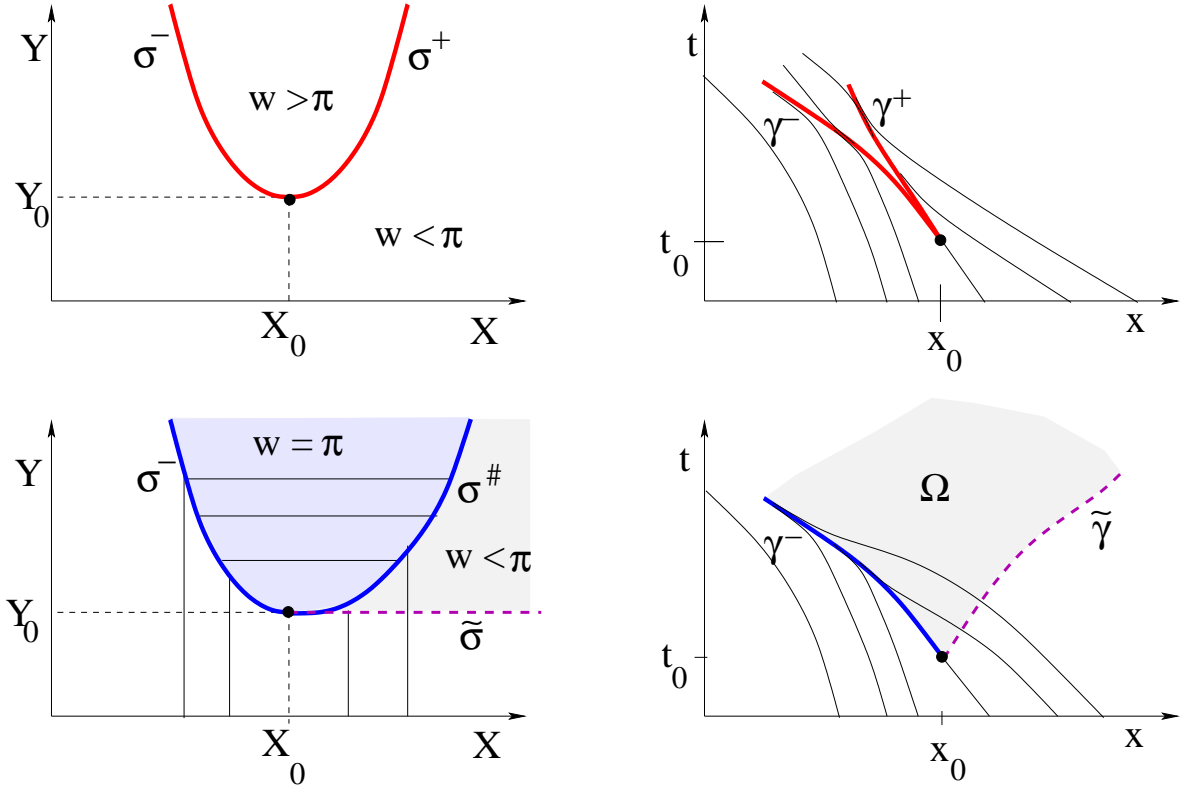}
\caption{ \small The positions of the singularities in the $X$-$Y$ plane
and in the $x$-$t$ plane. This refers to a point where a new singularity is formed, 
in the first family (i.e., for backward moving waves). 
Above: a conservative solution. Below: a dissipative solution.
Notice that the entire region between the curves 
$\sigma^-$ and $\sigma^\sharp$ 
is mapped onto the single curve $\gamma^-$. Indeed, horizontal segments
in the $X$-$Y$ plane are mapped into a single point.
In the $x$-$t$ plane, the two solutions differ only on the set
$\Omega$, bounded by the characteristic curves $\gamma^-$ 
(the image of 
both $\sigma^-$ and $\sigma^\sharp$) and $\tilde\gamma$ (the image of 
the line $\tilde\sigma$). }
\label{f:wa63}
\end{figure}

Up to the time $t_0= t(X_0, Y_0)$ where the  singularity appears, 
the conservative and the dissipative solution coincide.  
For $t>t_0$, they still coincide
outside the domain
\bel{Om}
\Omega~\doteq~\bigl\{(x,t)\,;\qquad t\geq t_0\,,\quad
\gamma^-(t)\leq x\leq \tilde\gamma(t)\bigr\},\eeq
where $\tilde \gamma$ is the forward characteristic through the point 
$(x_0, t_0)$.
Figure~\ref{f:wa63} shows the positions of 
these singularities in the $X$-$Y$ plane
and in the $x$-$t$ plane.
Figure~\ref{f:wa80} illustrates the difference in the profiles of 
the two solutions
for $t>t_0$.   
Our results  can be summarized as follows.
\v
{\bf Theorem 4.}
{\it In the above setting, 
the conservative solution $u^{cons}(t,\cdot)$ 
has two strong singularities at $x=\gamma^-(t)$ and $x=\gamma^+(t)$,
where $|u^{cons}_x|\to \infty$, and is smooth at all other points.

On the other hand, the dissipative solution $u^{diss}(t,\cdot)$ 
has a strong singularity at $x=\gamma^-(t)$,
where $|u^{diss}_x|\to \infty$, and a weak singularity 
along the forward characteristic $x=\tilde\gamma(t)$, where
$u^{diss}_x$ is continuous but the second derivative $u^{diss}_{xx}$
does not exist.

The difference between these two solutions can be estimated as
\bel{ucd}
\bigl\| u^{cons}(t,\cdot)-u^{diss}(t,\cdot)\|_{\C^0(\R)}
~=~\O(1)\cdot (t-t_0).\eeq
}
\v
\begin{figure}[htbp]
\centering
\includegraphics[width=0.8\textwidth]{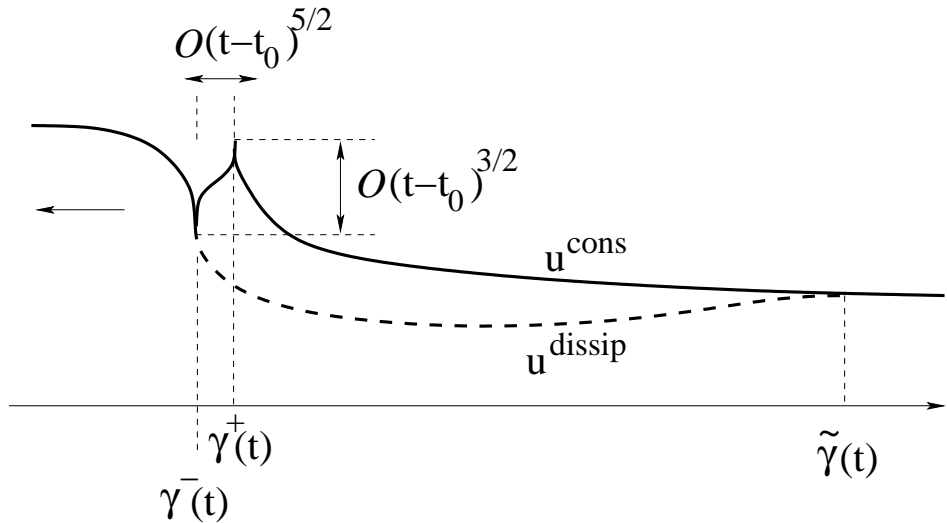}
\caption{ \small  Comparing a conservative and a dissipative solution,
at a time $t>t_0$, after a  singularity has appeared.
The  conservative solution has two strong singularities
at $\gamma^-(t)<\gamma^+(t)$,
while the dissipative solution has a strong singularity  
at $\gamma^-(t)$ and a weak singularity at $\tilde\gamma(t)$. 
The two solutions coincide for $x\leq \gamma^-(t)$ and for $x\geq
\tilde\gamma(t)$.}
\label{f:wa80}
\end{figure}

{\bf Proof.} {\bf 1.} To fix the ideas, 
assume that at the point $P=(X_0,Y_0)$
where the singularity is formed one has
$$w_{XX}~<~0,\qquad\qquad w_Y~>~0,\qquad\qquad c'(u)~>~0.$$
In the $X$-$Y$ coordinates, for smooth initial data 
the components $(x,t,u,w,z,p,q)$ of the conservative solution remain globally smooth.   On the other hand, for a dissipative solution
by (\ref{4.1})-(\ref{4.2}) we only know
that these components are  Lipschitz continuous.
\v
{\bf 2.} For $X\geq X_0$ we  denote by  $Y=\sigma^\sharp(X)$ the curve
where $w=\pi$, in the dissipative solution.
A Taylor approximation for $\sigma^\sharp$ is derived from
the identities
$$w(X,Y_0)~=~w_0 +  w_{XX}(X_0, Y_0)\cdot{(X-X_0)^2\over 2}
+\O(1)\cdot (X-X_0)^3,$$
$$w_Y(X,Y)~=~w_Y(X_0,Y_0)+\O(1)\cdot\bigl(|X-X_0|+|Y-Y_0|\bigr),$$
valid in the region where $w<\pi$.
Together, they imply
\bel{sse}
\sigma^\sharp(X)~=~Y_0 + 
\kappa(X-X_0)^2 + \O(1)\cdot (X-X_0)^3,\eeq
where $\kappa>0$ is the same constant found in (\ref{kdef})
for the conservative solution.

For $Y'>Y_0$, (\ref{wpi}) and (\ref{sse}) together imply
\bel{X-}
X^\sharp(Y')- X^-(Y')~=~2\left({Y'-Y_0\over \kappa}\right)^{1/2}+ \O(1) \cdot |Y'-Y_0|\,.
\eeq

\v

\begin{figure}[htbp]
\centering
\includegraphics[width=0.8\textwidth]{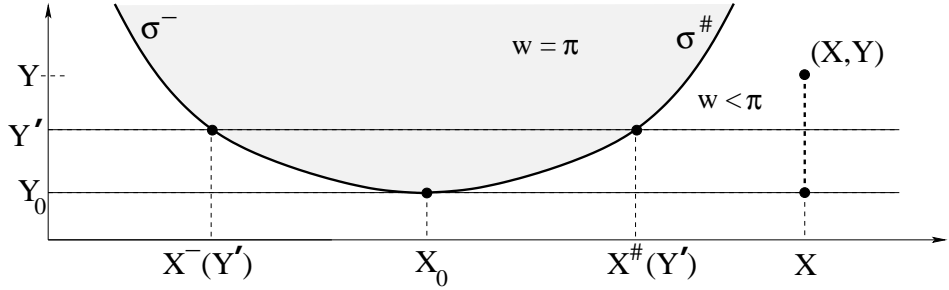}
\caption{ \small  Estimating the values of a dissipative solution
near a singularity. Notice that the functions
$x,t,u$ are constant on every horizontal segment 
contained in the shaded region where $w=\pi$.}
\label{f:wa81}
\end{figure}

{\bf 3.} Consider a point $(X,Y)$ with $X>X_0$ and 
$Y\leq \sigma^\sharp(X)$. By the second equation in (\ref{u}) it follows
\bel{u7}
u(X,Y)~=~u(X,Y_0) + \int_{Y_0}^Y \left( {\sin z\over 4c(u)}\, q
\right) (X,Y')\, dY'\,.\eeq
As  in Fig.~\ref{f:wa81}, for $Y'\in [Y_0, Y]$, 
call $X^-(Y')$ and $X^\sharp(Y')$
respectively  the points 
where $\sigma^-(X)=Y'$ and $\sigma^\sharp(X)=Y'$.
Since $z_X= q_X=0$ when $w=\pi$, by the second equations in
(\ref{4.1}) and in (\ref{4.2})
we have
\bel{z4}
\begin{split}
z(X, Y') ~&=~z(X^-(Y'), Y') + \int_{X^-(Y')}^{X} z_X(X',Y')\, dX'\\
&~=~z(X^-(Y'), Y')+
\int_{X^\sharp(Y')}^{X}\left({c'(u)\over 8c^2(u)}\,( \cos w - \cos z)
\,p\right)(X',Y')\, dX',
\end{split}
\eeq
\bel{q4}  q(X, Y') ~=~q(X^-(Y'), Y') + \int_{X^\sharp(Y')}^{X}
\left({c'(u)\over 8c^2(u)}\,
\big[\sin w-\sin z\big]\,pq\right)(X',Y')\, dX'.
\eeq
\v
{\bf 4.}
For notational convenience,  in the following we  denote by
$(x,t,u,w,z,p,q)(X,Y)$ the components describing a dissipative solution,
and by $(\hat x,\hat t,\hat u,\hat w,\hat z,\hat p,\hat q)(X,Y)$
the corresponding components of the conservative solution.
We observe that all these functions are Lipschitz continuous.
As shown in Fig.~\ref{f:wa63}, these two solutions can be different
only at points $(X,Y)$ in the region bounded by the curves $\sigma^-$ and $\tilde\sigma$, namely
$$\bigl\{ X\leq X_0\,,
\quad  Y>\sigma^-(X)\bigr\}~\cup~\bigl\{X\geq X_0\,,\quad Y>Y_0 \bigr\}.
$$
Consider a point $(X,Y')$ with $X>X_0$, $Y'<\sigma^\sharp(X)$.
By (\ref{4.1}), observing that $z=\hat z$ for $Y\leq Y_0$ and using (\ref{X-})
we find
\bel{z5}\bega{l}
\hat z(X,Y') - z(X, Y')\cr\cr
\quad \ds=~
\int_{X^-(Y')}^{X^\sharp(Y')} \hat z_X(X', Y')\, dX'
+\int_{X^\sharp(Y')}^{X}\left(\hat z_X-z_X\right)(X',Y')\, dX'
\cr\cr
\quad =~\ds
-{c'(u_0) (1+\cos z_0)\over 8c^2(u_0)} \, p_0\cdot 
\bigl( X^\sharp(Y')- X^-(Y')\bigr)
\cr\cr
\qquad\ds+\O(1)\cdot \bigl( X^\sharp(Y')- X^-(Y')\bigr)^2 +
 \O(1)\cdot  \bigl( X- X^\sharp(Y')\bigr)\, (Y'-Y_0)
 \cr\cr
\quad =~\ds
 -{c'(u_0) (1+\cos z_0)\over 4c^2(u_0)\,\kappa^{1/2}} \, p_0\cdot  (Y'-Y_0)^{1/2}
 +\O(1)\cdot |Y'-Y_0|\,.
\enda
\eeq
By (\ref{4.2}), a similar computation yields
\bel{q5}
\bega{l}
\hat q(X,Y') - q(X, Y')\cr\cr
\quad \ds=~
\int_{X^-(Y')}^{X^\sharp(Y')} \hat q_X(X', Y')\, dX'
+\int_{X^\sharp(Y')}^{X}\left(\hat q_X-q_X\right)(X',Y')\, dX'
\cr\cr
\quad =~\ds
-{c'(u_0) \sin z_0\over 8c^2(u_0)} \, p_0q_0\cdot 
\bigl( X^\sharp(Y')- X^-(Y')\bigr)
\cr\cr
\qquad\ds+\O(1)\cdot\bigl( X^\sharp(Y')- X^-(Y')\bigr)^2+ \O(1)\cdot  \bigl( X- X^\sharp(Y')\bigr)\, (Y'-Y_0)
\cr\cr
\quad =~\ds
-{c'(u_0) \sin z_0\over 4c^2(u_0)\, \kappa^{1/2}} \, p_0q_0\cdot  (Y'-Y_0)^{1/2}
 +\O(1)\cdot |Y'-Y_0|\,.
\enda
\eeq
Next, using the second equation in (\ref{u}) and recalling that 
$u(X,Y_0) = \hat u(X,Y_0)$, for any $X\in [X_0, \, X_0+1]$ 
and $Y\in [Y_0, \sigma^\sharp(X)]$ we obtain
\bel{u8}
\bega{l}
\ds\hat u(X,Y)-u(X,Y)
~ =~ \int_{Y_0}^Y \left({\sin \hat z\over 4c(\hat u)}\,\hat q
- {\sin z\over 4c(u)}\, q\right) (X,Y')\, dY'\cr\cr
\ds\qquad =~ \int_{Y_0}^Y \left({\cos z_0\over 4c(u_0)}\, 
q_0\cdot(\hat z-z)
+ {\sin z_0\over 4c(u_0)}\, \cdot(\hat q-q)
- {c'(u_0)\sin z_0\over 4c^2(u_0)}\, q_0\cdot(\hat u-u)\right)
 (X,Y')\, dY'\cr\cr
 \qquad\qquad\ds + \O(1)\cdot\int_{Y_0}^Y \Big(
 |\hat z-z|^2 + |\hat q-q|^2 + |\hat u-u|^2\Big)(X, Y')\, dY'\cr\cr
 \qquad\qquad\ds + \O(1)\cdot\Big( |X-X_0|+|Y-Y_0|\Big)
 \cdot\int_{Y_0}^Y \Big(
 |\hat z-z| + |\hat q-q| + |\hat u-u|\Big)(X, Y')\, dY'\cr\cr
 \qquad =~\eta_0\cdot (Y-Y_0)^{3/2} + \O(1)\cdot |Y-Y_0|^2,
\enda
\eeq
where the constant $\eta_0$ is computed by
\bel{eta0}
\eta_0~=~{2\over 3}\left[ - {c'(u_0) (1+\cos z_0)\over 4c^2(u_0) \,\kappa^{1/2}}\, p_0 \right]\cdot {\cos z_0\over 4c(u_0)}\, q_0 +
{2\over 3}\left[ - {c'(u_0) \,\sin z_0\over 4c^2(u_0) \,\kappa^{1/2}} \,p_0 q_0\right]\cdot 
{\sin z_0\over 4c(u_0)}\,.\eeq

{\bf 5.}  Using the second equations in (\ref{x}) and in (\ref{t}), we obtain similar estimates for the variables $x,t$.  Namely,
\bel{x8}
\bega{l}
\ds\hat x(X,Y)-x(X,Y)
~ =~ -\int_{Y_0}^Y \left({1+\cos \hat z\over 4}\,\hat q
- {1+\cos z\over 4}\, q\right) (X,Y')\, dY'\cr\cr
\ds\qquad 
=~ \int_{Y_0}^Y \left({\sin z_0\over 4}\, 
q_0\cdot(\hat z-z)
- {1+\cos z_0\over 4}\, \cdot(\hat q-q)
\right)
 (X,Y')\, dY'\cr\cr
 \qquad\qquad\ds + \O(1)\cdot\int_{Y_0}^Y \Big(
 |\hat z-z|^2 + |\hat q-q|^2 \Big)(X, Y')\, dY'\cr\cr
 \qquad\qquad\ds + \O(1)\cdot\Big( |X-X_0|+|Y-Y_0|\Big)
 \cdot\int_{Y_0}^Y \Big(
 |\hat z-z| + |\hat q-q|\Big)(X, Y')\, dY'\cr\cr
 \qquad =~ \O(1)\cdot |Y-Y_0|^2.
\enda
\eeq
Indeed, the coefficient of the leading order term $\O(1)\cdot (Y-Y_0)^{3/2}$ 
vanishes.  Similarly, 
\bel{t8}
\bega{l}
\ds\hat t(X,Y)-t(X,Y)
~ =~ -\int_{Y_0}^Y \left({1+\cos \hat z\over 4c(\hat u)}\,\hat q
- {1+\cos z\over 4c(u)}\, q\right) (X,Y')\, dY'\cr\cr
\ds
=~ \int_{Y_0}^Y \left({\sin z_0\over 4c(u_0)}\, 
q_0\cdot(\hat z-z)
- {1+\cos z_0\over 4c(u_0)}\, \cdot(\hat q-q) + {1+\cos z_0\over 4c^2(u_0)}c'(u_0)
q_0\cdot (\hat u - u)
\right)
 (X,Y')\, dY'\cr\cr
 \qquad\qquad\ds + \O(1)\cdot\int_{Y_0}^Y \Big(
 |\hat z-z|^2 + |\hat q-q|^2+ |\hat u - u|^2 \Big)(X, Y')\, dY'\cr\cr
 \qquad\qquad\ds + \O(1)\cdot\Big( |X-X_0|+|Y-Y_0|\Big)
 \cdot\int_{Y_0}^Y \Big(
 |\hat z-z| + |\hat q-q|+|\hat u -u|\Big)(X, Y')\, dY'\cr\cr
=~ \O(1)\cdot |Y-Y_0|^2.
\enda
\eeq
\v
%
%
{\bf 6.} The  estimate (\ref{u8}) provides a 
bound on the difference $\hat u - u$ 
between  a conservative and a dissipative solution, at a given point
$(X,Y)$.
However, our main goal is to estimate the difference $\hat u - u$ as functions
of the original variables $x,t$.   For this purpose, consider a dissipative solution
$u$ and a point 
\bel{P} P~=~(x, t)~=~ \bigl(x(X, Y),\, t(X, Y)\bigr),\eeq with 
\bel{Omega}X~>~X_0, \qquad\qquad Y_0~<~Y~<~\sigma^\sharp(X).\eeq
Moreover, let $\hat u$ be the conservative solution with the 
same initial data, and
let $(\Hat X, \Hat Y)$ be the point 
which is mapped to $P$ in the conservative solution, so that
\bel{htx}
P~=~(x, t)~=~ \bigl(\hat x(\Hat X, \,\Hat Y), 
\hat t(\Hat X, \Hat Y)\bigr).\eeq
Using (\ref{u8}), (\ref{x8}), (\ref{t8}),  and recalling 
that the conservative solution $\hat u= \hat u(x,t)$ is H\"older continuous
of exponent $1/2$ w.r.t.~both variables $x,t$, we obtain
\bel{uhu}\bega{l}
|\hat u(x,t) - u(x,t)|~\leq~\bigl|\hat u(\Hat X,\Hat Y)- \hat u(X,Y)
\bigr| + 
\bigl|\hat u(X,  Y)- u(X, Y)\bigr|\\[4mm]
\qquad =~\O(1)\cdot \Big(\bigl|\hat x(\Hat X,\Hat Y)- \hat x(X,Y)
\bigr|^{1/2}+\bigl|\hat t(\Hat X,\Hat Y)- \hat t(X,Y)
\bigr|^{1/2}\Big) + \O(1) \cdot |Y-Y_0|^{3/2}
\\[4mm]
\qquad =~\O(1)\cdot \Big(\bigl| x(X,Y)- \hat x(X,Y)
\bigr|^{1/2}+\bigl|t(X, Y)- \hat t(X,Y)
\bigr|^{1/2}\Big) + \O(1) \cdot |Y-Y_0|^{3/2}
\\[4mm]
\qquad =~\O(1)\cdot |Y-Y_0|
.\enda\eeq
In a neighborhood of $(x_0, t_0)$
we have
\bel{tY}t_Y ~=~{1+\cos z\over 4c(u)}\,q ~>~ {1+\cos z_0\over 5 c(u_0)}q_0~
>0\,.\eeq
For $(x,t)$ as in (\ref{P})-(\ref{Omega}), one has
$$t-t_0~=~\bigl[t(X,Y) - t(X, Y_0)\bigr] +
 \bigl[t(X, Y_0) - t(X_0, Y_0)\bigr]~\geq~
{1+\cos z_0\over 5 c(u_0)}q_0\cdot |Y-Y_0|\,.$$
Together with (\ref{uhu}), this proves (\ref{ucd}).
\v
{\bf 7.} It remains to prove that the solution $u=u(t,x)$ is not twice
differentiable along the forward characteristic $\tilde\gamma$.

Consider a point $(x_1, t_1) = \bigl(x(X,Y_0), t(X, Y_0)\bigr)$
on $\tilde \gamma$, with $X>X_0$.  
Let $x=\gamma(t)$ be the backward characteristic through
$(x_1, t_1)$, so that 
$$\gamma(t_1) = x_1,\qquad\qquad \dot\gamma(t)~=~c(u(\gamma(t),t)).$$
Assume that $u$ were twice differentiable at the point $(x_1, t_1)$.
Then the map $t\mapsto u(\gamma(t), t)$ would also 
be twice differentiable at $t=t_1$. Indeed
$${d\over dt} u\bigl(\gamma(t), t\bigr) ~=~- c(u)u_x + u_t\,,$$ 
\bel{ddt}{d^2\over dt^2} u\bigl(\gamma(t), t\bigr) ~
=~- c'(u) \, \bigl(-c(u) u_x+u_t\bigr)u_x 
+ c^2(u)u_{xx} -2c(u) u_{xt} + u_{tt}\,.\eeq
To reach a contradiction, 
consider the map $\tau\mapsto Y(\tau)$ 
implicitly defined by 
$$t(X, Y(\tau))~=~\tau\,.$$
By (\ref{tY})  this map is well defined.  In particular, $Y(t_1) = Y_0$.
We thus have
\bel{dt}{d\over dt} u\bigl(\gamma(t), t\bigr) ~=~{d\over dt} u(X, Y(t))
~=~u_Y\cdot {4c(u)\over (1+\cos z)q}~=~{\sin z\over 1+\cos z}\,.\eeq
We now show that this first derivative cannot be a Lipschitz continuous
function of time, for $t\approx t_1$.   Indeed, by (\ref{dt}) 
and the mean value theorem we have
\bel{ddd}\bega{rl}\ds{d\over dt} u(X, Y(t)) - {d\over dt} u(X, Y(t_1))&=~
\ds{\sin z(X, Y(t))\over 1+\cos z(X, Y(t))}- 
{\sin z(X, Y_0)\over 1+\cos z(X, Y_0)}\cr\cr
&=~\ds {1\over 1+\cos z(X, Y^\sharp)}\cdot \bigl[ z(X,Y(t)) - z(X,
Y_0)\bigr],
\enda
\eeq
for some intermediate value $Y^\sharp\in [Y_0, Y(t)]$.
 Call $\hat z= \hat z(X,Y)$ the corresponding conservative solution.
 Observe that $\hat z$ is smooth and coincides with $z$ 
 on the horizontal
 line $\{ Y=Y_0\}$.   
Using (\ref{z5}) we obtain
 $$\bega{rl} \bigl|z(X,Y(t)) - z(X, Y_0)
 \bigr|&\geq ~\bigl[z(X,Y(t))- \hat z(X, Y(t))
 \bigr] -\bigl[\hat z(X,Y(t))- \hat z(X, Y_0)
 \bigr]\cr\cr
 &\geq~\ds
  {c'(u_0) (1+\cos z_0)\, p_0\over 4c^2(u_0)\,\kappa^{1/2}} 
  \cdot  (Y(t)-Y_0)^{1/2}
 -\O(1)\cdot |Y(t)-Y_0|\,.\enda
 $$
As a consequence, for  $t\approx t_1$, 
the function $t\mapsto z(X, Y(t))$ is not Lipschitz continuous,
and the same applies to the left hand side of (\ref{ddd}).
We thus conclude that the map  $t\mapsto u(\gamma(t),t)$ cannot be twice
differentiable at $t=t_1$, in contradiction with (\ref{ddt}).
This completes the proof of Theorem~4.
\endproof

 \v
{\bf Acknowledgment.} This research  was partially supported
by NSF, with grant  DMS-1411786: ``Hyperbolic Conservation Laws and Applications".

\v

\end{document}